\documentclass[11pt,oneside]{article}
\usepackage{amssymb,amsmath}
\usepackage{url}
\renewcommand{\le}{\leqslant}
\renewcommand{\ge}{\geqslant}

\newcommand{\mad}{\mathbf{Mad}}

\newcommand{\RR}{\mathbb{R}}
\newcommand{\ZZ}{\mathbb{Z}}
\newcommand{\QQ}{\mathbb{Q}}

\newcommand{\NN}{\mathbb{N}}

\renewcommand{\epsilon}{\varepsilon}

\newtheorem{lemma}{Lemma}
\newtheorem{theorem}{Theorem}
\newtheorem{proposition}{Proposition}

\newenvironment{proof}{\textsc{Proof. }}{\ \newline\hspace*{\fill}$\boxtimes$}

\begin{document}
\title{Estimating lower limit in the $p$-adic Littlewood conjecture}

\author{Dmitry Badziahin}


\maketitle

\begin{abstract}
We verify that $\liminf_{q\to\infty} q\cdot |q|_p\cdot
||qx||<\epsilon$ for all real $x$, small primes $p$ and relatively
small $\epsilon$. This result supports the famous $p$-adic
Littlewood conjecture which states that the above lower limit is
equal to 0 for all $x\in\RR$. In particular, the result is
established for $p=2$ with $\epsilon=1/25$. For $3\le p\le 29$, the
upper bounds for $\epsilon$ vary, but they are always at most
$1/10$.
\end{abstract}

{\footnotesize{{\em Keywords}: $p$-adic Littlewood conjecture

Math Subject Classification 2020: 11J13}}

\section{Introduction}

The $p$-adic Littlewood conjecture (PLC) was introduced in 2004 by
de~Mathan and Teuli\'e~\cite{mat_teu_2004}. It states that for all
primes $p$ and real numbers $x$ one has
$$
\liminf_{q\to\infty} q\cdot |q|_p\cdot ||qx||=0,
$$
where $||\cdot||$ denotes the distance to the nearest integer.
Despite the efforts of many mathematicians, the conjecture remains
open. We refer the reader to~\cite{ein_kle_2007, bbek_2015,
bkn_2024} for recent progress regarding the PLC. In particular, it
is known that for all values of $p$ the set of counterexamples $x$
has zero Hausdorff dimension and every such value $x$ must satisfy
several restrictive conditions.

In~\cite{badziahin_2016}, the author developed an algorithm which,
for a given $\epsilon>0$, verifies the inequality
\begin{equation}\label{eq1}
\liminf_{q\to\infty} q\cdot |q|_p \cdot ||qx||<\epsilon
\end{equation}
for all $x\in\RR$. Essentially, for every small interval with
endpoints in consecutive Farey fractions the algorithm checks
whether all its points satisfy~\eqref{eq1}. Such an algorithm
returns a positive answer if~\eqref{eq1} is satisfied for all
$x\in\RR$. On the other hand, if there exists a counterexample
$x\in\RR$ that violates the inequality, the algorithm would never
terminate.

Applied to the 2-adic LC, the algorithm provided the following
results \cite[Table 2]{badziahin_2016}:
\begin{center}
\begin{tabular}{|c|c|c|c|}
\multicolumn{4}{c}{{\bf Table 1.} Output for 2-adic LC from~\cite{badziahin_2016}.}\\
\hline $\epsilon$&$Q$&Number of iterations&Time
spent\\
\hline 1/4&8&6&0\\
\hline 1/5&16&11&0\\
\hline 1/6&128&18&0\\
\hline 1/7&1024&23&0\\
\hline 1/8&10240&34&0.234 s\\
\hline 1/9&65536&39&3.02 s\\
\hline 1/10&$\ge$446431232&$>$54&$>$60 h\\
\hline
\end{tabular}
\end{center}
Here, $Q = Q(\epsilon)\in \NN$ is defined as the smallest integer
such that $\min_{1\le q\le Q} \{q\cdot |q|_p\cdot ||qx||\}<\epsilon$
for all $x\in \RR$. Notice a significant increase in the values of
$Q$ when $\epsilon$ decreases from $1/9$ to $1/10$. This phenomenon,
together with the recent disproof of an analogous $t$-adic
Littlewood conjecture over finite function fields~\cite{aln_2021},
raises doubts about the validity of the PLC as well. To gather more
evidence for or against the conjecture, we continue the search and
verify the $p$-adic Littlewood conjecture for $p=2$ and
substantially smaller values of $\epsilon$. The main outcome of this
search is the following result:

\begin{theorem}\label{th1}
For all $x\in\RR$ one has
$$
\liminf_{q\to\infty} q\cdot |q|_2\cdot ||qx|| <1/25.
$$
\end{theorem}

We also verify the conjecture for all primes $p\le 29$ and various
values of $\epsilon$. The exact values of $\epsilon$ for each prime
$p$ are provided in the Numerical Results section. However, in all
cases one has $\epsilon\le 1/10$.

Of course, no matter how small the value of $\epsilon>0$ we can use
in Theorem~\ref{th1}, this still does not constitute the proof of
the PLC. However, the results of this search lend further support to
the conjecture: not only do we reduce $\epsilon$ by a factor greater
than 2, but we also no longer observe the large jumps seen in
Table~1. To achieve such an improvement, we adopt a new approach
that leverages the dynamical nature of the problem. We describe it
in detail in the next section.

{\bf Remark.} While completing this note, we discover a similar work
on this topic~\cite{vit_vuk_2025} where Theorem~\ref{th1} was
verified for $\epsilon=1/15$. We believe it is worth mentioning
here.

\section{Description of the algorithm}

First of all, instead of considering all real values $x$, we look at
them as elements of $\RR/\ZZ$. Then we make an observation: if $x$
is a counterexample to~\eqref{eq1} for a given $\epsilon$ then so is
the value $px$ for the same $\epsilon$. In other words, the set of
counterexamples to the PLC with the same parameter $\epsilon$ is
invariant under multiplication by $p$.

Next, we establish an auxiliary statement.

\begin{lemma}\label{lem3}
Given $\epsilon>0$, if there exists $x\in\RR/\ZZ$ that satisfies
$$
\liminf_{q\to\infty} q\cdot |q|_p \cdot ||qx||\ge \epsilon
$$
then there exists $y\in\RR/\ZZ$ that satisfies
\begin{equation}\label{lem3_eq}
\inf_{q\in\NN}q\cdot |q|_p \cdot ||qy||\ge \epsilon.
\end{equation}
\end{lemma}

\proof Consider $x$ that satisfies the condition of the lemma and
consider any limit point $y$ of its orbit $O(x):= \{p^nx :
n\in\NN\}\subset \RR/\ZZ$. Suppose that
$$
q\cdot |q|_p \cdot ||qy||= \epsilon-\delta<\epsilon
$$
for some $q\in\NN$. Then by construction, there exists an infinite
sequence of indices
$$
k_1<k_2<\ldots <k_i <\ldots
$$
such that $p^{k_i} x$ is so close to $y$ that
$$
q\cdot |q|_p\cdot ||q\cdot p^{k_i}x|| \le
\epsilon-\delta/2<\epsilon,\quad i\in\NN.
$$
Therefore
$$
p^{k_i}q\cdot |p^{k_i}q|_p\cdot ||p^{k_i}qx||\le \epsilon-\delta/2
$$
and finally
$$
\liminf_{q\to\infty} q\cdot |q|_p ||qx|| <\epsilon,
$$
which is a contradiction.
\endproof

Due to Lemma~\ref{lem3}, it is sufficient to check if the following
set is nonemtpy:
$$
\mad_p(\epsilon):= \{x\in \RR/\ZZ: \inf_{q\in \NN} q\cdot |q|_p\cdot
||qx||\ge \epsilon\}.
$$
Clearly, it is also invariant under multiplication by $p$. Moreover,
as a countable intersection of closed sets, it is closed.

Consider $x\in\mad_p(\epsilon)$ and look at its orbit $O(x)$. Let
$f:\RR/\ZZ\to \RR^+$ be a continuous function. We call $x$ an
$f$-bottom number if for all $n\ge 1$, $f(p^nx)\ge f(x)$. Notice
that $\overline{O(x)}$ contains at least one $f$-bottom number.
Since $x\in\mad_p(\epsilon)$ implies $\overline{O(x)}\subseteq
\mad_p(\epsilon)$, we derive the following
\begin{proposition}\label{prop2}
For any continuous function $f:\RR/\ZZ \to \RR^+$, if
$\mad_p(\epsilon)$ is non-empty then there exists an $f$-bottom
number $y\in \mad_p(\epsilon)$.
\end{proposition}

Now, one of the key ideas of the algorithm is to verify that the set
$\mad_p(\epsilon)$ does not contain $f$-bottom numbers for a
suitably chosen continuous function $f$. We always use $f(z) =
||z||$, however a smarter choice of $f$ may potentially make an
algorithm faster.

Another observation is that, due to symmetry, $x\in\mad_p(\epsilon)$
if and only if $-x \in\mad_p(\epsilon)$. The function $f(z) = ||z||$
has the same symmetry as well. Therefore it is sufficient to verify
that $\mad_p(\epsilon)$ has no $f$-bottom numbers in the interval
$(0,1/2)$.

\begin{lemma}\label{lem1}
Let $p$ be prime and $f(x) = ||x||$. If
\begin{equation}\label{lem1_eq}
x\in(0,1/2)\cap \left(\frac{c}{p^n+1},\frac{c}{p^n-1}\right)
\end{equation}
for some $c,n\in\NN$ then $x$ is not an $f$-bottom number.
\end{lemma}

\proof Indeed, if $x$ satisfies~\eqref{lem1_eq} then $||p^nx|| =
|p^nx - c| < x = ||x||$ and hence $x$ is not an $f$-bottom number.
\endproof

We say that intervals $\left(\frac{c}{p^n+1},\frac{c}{p^n-1}\right)$
from~\eqref{lem1_eq} are of type~1 and denote them by $I_1(c,n)$.

\begin{lemma}\label{lem2}
Let $p$ be prime. If
\begin{equation}\label{lem2_eq}
x\in \left(\frac{c}{p^n d} - \frac{\epsilon}{p^n d^2},
\frac{c}{p^nd} + \frac{\epsilon}{p^n d^2}\right) =: I_2(c,d,n) =
I_2^\epsilon(c,d,n)
\end{equation}
for some $c,d,n\in\NN$ then
$$
\inf_{q\in \NN} q\cdot |q|_p\cdot ||qx||< \epsilon.
$$
\end{lemma}

Lemma~\ref{lem2} follows immediately from the inequality $p^nd\cdot
|p^nd|_p\cdot ||p^ndx||<\epsilon$ for all $x$ in the
interval~\eqref{lem2_eq}. We say that all such intervals
$I_2(c,d,n)$ are of type~2. To simplify the implementation of the
algorithm we always set $\epsilon = 1/E$ where $E$ is a positive
integer. Then $I_2(c,d,n)$ can be rewritten as
$$
I_2(c,d,n) =
\left(\frac{cdE-1}{Ep^nd^2},\frac{cdE+1}{Ep^nd^2}\right).
$$

Since rational values of $x$ obviously satisfy PLC, for
$\epsilon\in\QQ$ we can slightly strengthen Lemmata~\ref{lem1}
and~\ref{lem2} by considering closures of both intervals
in~\eqref{lem1_eq} and~\eqref{lem2_eq}.

\begin{proposition}\label{prop1}
Fix $\epsilon>0$. Suppose that the interval $[0,1/2]$ is covered by
a union of intervals of type~1 and~2. Then~\eqref{eq1} is satisfied
for all $x\in\RR/\ZZ$.
\end{proposition}

\proof By Lemma~\ref{lem3}, we need to check that $\mad_p(\epsilon)$
is empty. Suppose it is not. Then choose an $f$-bottom number $y\in
\mad_p(\epsilon)$. Due to Proposition~\ref{prop2}, such a number
exists. If $y>1/2$, replace it with $1-y$. By Lemma~\ref{lem1}, $y$
is not covered by any interval of type~1. Hence it is covered by
some interval $I_2(c,d,n)$ and then Lemma~\ref{lem2} implies that
$$
p^nd\cdot |p^nd|_p\cdot ||p^ndy|| <\epsilon.
$$
This contradicts the choice of $y\in\mad_p(\epsilon)$.
\endproof

To prove the PLC for a given $\epsilon$, we build a cover of $[0,
1/2]$ by intervals of type~1 and~2 in the following way.
\begin{enumerate}
\item Start with $a/b = 1/2$.
\item Consider a number $a/b$ and look for an interval $I$ of type
1 or 2 such that $a/b\in I$ and $a/b$ is not the left endpoint of
$I$.
\item If one finds several such intervals $I$, keep one with the
largest length $|I|$.
\item Set $a/b$ to be the left endpoint of $I$.
\item If $a/b>0$ go back to step~2 and repeat the iteration.
\end{enumerate}

Verification of wether $a/b$ belongs to an interval $I$ of type~1
can be done in the following way. Notice that $p^n \frac{a}{b} =
\frac{a_n}{b}$, where $a_n \equiv p^na$ (mod $b$) i.e. $a_n$ equals
the remainder of $p^na$ after division by $b$. Then the condition
$$\left|\left|p^n \frac{a}{b}\right|\right| \le \left|\left|\frac{a}{b}\right|\right|$$ is satisfied if $a_n<a$ or $b-a_n < a$. Hence, schematically the
test will look as follows

\begin{itemize}
\item[]{\bf Step 1.} Set $a_n:=a$
\item[]{\bf Step 2.} For $n$ from 1 to some fixed upper bound
$\mathbf{Max_n}$ do the cycle

\begin{itemize}
\item[]{\bf Step 3.} Set $a_n:= p a_n$ and reduce it modulo $b$
\item[]{\bf Step 4.} If $a_n<a$ or $b-a_n \le a$ then we have found
$I$ of type~1 that covers $a/b$. In that case we perform the
following operations
\begin{itemize}
\item If $a_n<a$ set $Left\_a:= \lfloor p^n a/b\rfloor$ else
$Left\_a:=\lfloor p^na/b\rfloor + 1$. Here we denote by
$Left\_a/Left\_b$ the left end point of $I$.
\item Set $Left\_b:= p^n + 1$
\item Set $S := p^n$. This variable will indicate the size of $I$.
It approximately equals $1/|I|$. To save time, we do not compute its
precise value but use some rough approximation. Indeed,
$$
\frac{c}{p^n-1} - \frac{c}{p^n+1} = \frac{2c}{(p^n-1)(p^n+1)} \asymp
\frac{1}{p^n}
$$
if we assume that $a/b\asymp 1$.
\item Leave the cycle
\end{itemize}

\end{itemize}

\end{itemize}

Now, focus on the intervals of type~2. Notice that by Legendre's
theorem, if $\epsilon<\frac12$ and $a/b\in I_2(c,d,n)$ then $c/d$ is
a convergent of $\frac{p^na}{b}$. Therefore it is sufficient to
restrict the search to such numbers $c/d$. Also note that
$|I_2(c,d,n)| = \frac{2\epsilon}{p^nd^2}$ therefore we can restrict
our search to those intervals $I_2(c,d,n)$ that satisfy $p^nd^2\le
\epsilon S$. Schematically, the verification can be done in the
following way.


\begin{itemize}
\item[]{\bf Step 1.} Go through all integer values $n$, starting
from $n=0$, until $p^n$ becomes larger than $\epsilon S$, increasing
$n$ by one after each iteration. The iteration is

\begin{itemize}
\item[]{\bf Step 2.} Go through all convergents $c/d$ of $\frac{p^na}{b}$
until $p^nd^2$ exceeds $\epsilon S$. Each iteration consists of the
following check:

\begin{itemize}
\item[]{\bf Step 3.} If $d\cdot |p^nad - bc| < \epsilon b$ then do
the following

\begin{itemize}
\item Factor out all primes $p$ from $d$. That is, represent $p^nd =
p^{n^*}d^*$ where $n^*,d^*\in\ZZ$ and $p\nmid d^*$.
\item Set $Left\_a:= cd^*E-1$ and $Left\_b:=
Ep^{n^*}(d^*)^2$.
\item Set $S:=Ep^{n^*}(d^*)^2$.
\end{itemize}

\end{itemize}
\end{itemize}
\end{itemize}

After we complete all the checks, we update $a:=Left\_a$ and
$b:=Left\_b$.

\section{Small example}

Before making $\epsilon$ as small as possible, let's run the
algorithm for relatively large $\epsilon = \frac18$ and $p=2$. Such
a small example will help us better understand how the algorithm
works.

We start with $a/b=1/2$. We start with the search for intervals of
type~1 that cover $a/b$. By a quick check
$$
2a \;\mathrm{mod}\; b = 2\; \mathrm{mod}\; 2 = 0 < a
$$
infers that $1/2\in I_1(1,1) = \big[\frac13, 1\big]$. Set $S := 2$.

Since $1 = 2^0$ is already bigger than $\epsilon S = 2/8$, we skip
the search for intervals of type~2 and choose $I_1(1,1)$  as the
first interval in the cover.

Next, we let $a/b = 1/3$, the left endpoint of $I_1(1,1)$, and
search for a covering interval of this number. Similar checks as
above show that $1/3\in I_1(1,2) = \big[\frac15, \frac13\big]$. Then
$S = 4$ and we still have $1> \epsilon S$. Therefore we skip the
search for intervals of type~2.

In the next iteration, $a/b=1/5$. We verify that $1/5\in I_1(3,4) =
\big[\frac{3}{17}, \frac{3}{15}\big]$.

Continuing like that with $a/b=3/17$, we observe that this number is
 in $I_1(45,8) = \big[\frac{45}{257}, \frac{45}{255}\big]$. This
time $S = 256$ and we need to search for intervals $I_2(c,d,n)$ of
type~2 that cover $a/b$. However, according to Step~2 of the search,
their parameters $d$ and $n$ should satisfy $2^nd^2\le \epsilon S =
\frac{256}{8} = 32$. This in turn implies that $n\le 5$ and $d\le 5$
and a very quick search confirms that there are no suitable
intervals of type~2 that cover $3/17$.

Now we have $a/b = \frac{45}{257}$. We observe that $\frac{45}{257}
\in I_1(11475, 16)$, set $S = 2^{16} = 65536$ and continue with the
search of intervals of type~2. One gets that
$$\frac{45}{257} = [0;5,1,2,2,6]$$
and checks that there are no appropriate convergents $c/d$ close
enough to $a/b$ which satisfy $d^2 < \epsilon 2^{16}$. Next, we
compute
$$\frac{2\cdot 45}{257} = [0;2,1,5,1,12]$$
and verify that the convergent $c/d = [0;2,1,5,1] = 7/20$ is close
enough to $2a/b$. By factoring out all divisors 2 from the
denominator 20 of the convergent, we end up with $c/d = 7/5$ which
is close enough to $2^3a/b$ and therefore $\frac{45}{257} \in
I_2(7,5,3) = \big[\frac{279}{1600}, \frac{281}{1600}\big]$. The
bound $S$ after that reduces from $2^{16}$ to $1600$. Further search
of parameters $n\ge 2$ and convergents $c/d$ does not find any
better interval.

Continuing as above, we end up with the following cover of
$[0,1/2]$:
$$
\begin{array}{rl}
(0,1/2)\in& I_1(1,1)\cup I_1(1,2)\cup I_1(3,4)\cup I_1(45,8)\cup
I_2(7,5,3)\cup I_1(11,6)\cup I_2(1,3,1)\cup I_1(5,5)\cup \\
&I_1(155,10)\cup I_2(3,5,2)\cup I_1(19,7)\cup I_2(33,7,5)\cup
I_1(603,12)\cup I_2(5,17,1)\cup\\
&I_1(1203,13)\cup I_2(47,5,6)\cup I_1(2405,14)\cup I_2(16,109,0)\cup
I_1(38477,18)\cup\\
&I_2(155,33,5)\cup I_2(1503,5,13)\cup I_1(9619,16)\cup
I_2(263,7,8)\cup I_1(75,9)\cup I_2(7,3,4)\cup\\
&I_1(37,8)\cup I_2(1,7,0)\cup I_1(1,3)\cup I_2(0,1,0).
\end{array}
$$
In total, we have that the interval $[0,1/2]$ is covered by 29
intervals, 16 of them are of type~1 and the remaining 13 are of
type~2. Finally, Proposition~\ref{prop1} then implies~\eqref{eq1}
for all real $x$ and $\epsilon=\frac18$. Notice that intervals of
type~1 do not depend on $\epsilon$. Therefore they will always cover
a substantial part of $[0,1/2]$ no matter how small $\epsilon$ we
choose.

\section{Numerical results}

The algorithm was initially implemented on C++ with use of NTL
library that handles arbitrary large integers. For $p=2$ it was
launched on one core of AMD Ryzen 9 5990x CPU with parameters
$\epsilon=1/E$ where $E\in \NN$, $10\le E\le 23$. The algorithm for
$p=2$ was then rewritten with NTL being replaced by GMP library.
That made it approximately 3 times faster. A modified code was
launched on the same CPU for $\epsilon = 1/24$. Finally, for
$\epsilon=1/25$ and $p=2$ the interval $[0,1/2]$ was divided into 12
pieces and a cover for each piece was built in parallel on 12 cores
of the same CPU. Because of the slight cover overlap, we can not
provide the exact number of intervals in it but instead provide the
approximate number. The results of the algorithm are provided in
Table~2 below.

The calculations suggest that the number of intervals in the cover
grow exponentially with $E$. It does not show any irregularities as
in Table~1.

\begin{center}
\begin{tabular}{|c|c|c|}
\multicolumn{3}{c}{{\bf Table 2.} Results for 2-adic LC.}\\
\hline $\epsilon$&Number of intervals in the cover&Time spent\\
\hline 1/10&197&0 s\\
\hline 1/11&421&0 s\\
\hline 1/12&971&0 s\\
\hline 1/13&2726&0 s\\
\hline 1/14&10407&2.25 s\\
\hline 1/15&34847&7.84 s\\
\hline 1/16&84424&21.77 s\\
\hline 1/17&251781&1 m 9 s\\
\hline 1/21&44199404&6 h 4 m\\
\hline 1/22&124564436&19 h 53 m\\
\hline 1/23&478809708&84 h 11 m\\
\hline 1/24&1988129786&82 h 54 m\\
\hline 1/25&$\approx$7369690680&44 h 41m\\
\hline
\end{tabular}
\end{center}

Next, we run the algorithm for various small primes and values of
$\epsilon$. This time, we did not aim to minimise $\epsilon$ or
optimise the code. Our primary goal was to rule out potential
counterexamples to the PLC for $p>2$ and moderately small values of
$\epsilon$ where any such counterexample would cause the algorithm
to quickly stall. The results of the search are presented in
Table~3. For each $p$ and each value $\epsilon$ tested, the
algorithm successfully verified that the $PLC$ holds with that
$\epsilon$. These results further support the conjecture for all
primes~$p$.

\begin{center}
\begin{tabular}{|c|c|c|c|}
\multicolumn{4}{c}{{\bf Table 3.} Results for general $p$-adic LC.}\\
\hline $p$&$\epsilon$&Number of intervals&Time spent\\
\hline 3&1/18&3526231&13 m 29 s\\
\hline 5&1/15&837078&2 m 1 s\\
\hline 7&1/15&4082365&9 m 6 s\\
\hline 11&1/13&3871350&7 m 53 s\\
\hline 13&1/12&1328134&2 m 27 s\\
\hline 17&1/12&3809979&6 m 54 s\\
\hline 19&1/12&8579511&15 m 19 s\\
\hline 23&1/11&3674595&5 m 28 s\\
\hline 29&1/10&4805531&8 m 19 s\\
\hline
\end{tabular}
\end{center}

\bigskip
\noindent Dzmitry Badziahin\\ \noindent The University of Sydney\\
\noindent Camperdown 2006, NSW (Australia)\\
\noindent {\tt dzmitry.badziahin@sydney.edu.au}

\end{document}